\documentclass{amsart}

\usepackage{amssymb}
\usepackage{mathtools}
\usepackage{microtype}
\usepackage{xcolor}
\usepackage[hidelinks]{hyperref}

\theoremstyle{plain}
\newtheorem{theorem}{Theorem}[section]
\newtheorem{lemma}[theorem]{Lemma}
\newtheorem{corollary}[theorem]{Corollary}
\theoremstyle{remark}

\numberwithin{equation}{section}

\newcommand{\defeq}{\mathrel{:=}}
\newcommand{\E}{\mathbb E}
\newcommand{\PP}{\mathbb P}
\newcommand{\Var}{\operatorname{Var}}
\newcommand{\TV}[2]{\mathrm d_{\mathrm{TV}}\!\left(#1,#2\right)}
\newcommand{\tmix}{t_{\mathrm{mix}}^{(n)}}
\newcommand{\cL}{\mathcal L}
\newcommand{\zero}{\mathbf 0}
\newcommand{\one}{\mathbf 1}

\title[Mixing and cutoff for the noisy SIS model]
      {Mixing time and cutoff for the noisy SIS model on graphs}

\author[W. R. KhudaBukhsh]{Wasiur R. KhudaBukhsh}
\address{School of Mathematical Sciences, University of Nottingham,
Nottingham, United Kingdom}
\email{wasiur.khudabukhsh@nottingham.ac.uk}

\author[Y. Xiang]{Yangrui Xiang}
\address{Department of Mathematics, Louisiana State University,
Baton Rouge, Louisiana, USA}
\email{yxiang8@lsu.edu}
\thanks{Yangrui Xiang is the corresponding author.}

\date{}

\hypersetup{
  pdfauthor={Wasiur R. KhudaBukhsh and Yangrui Xiang},
  pdftitle={Mixing time and cutoff for the noisy SIS model on graphs}}

\begin{document}

\begin{abstract}
We study the mixing time and cutoff phenomenon for the con\-tin\-u\-ous-time susceptible--infected--susceptible (SIS) model with external sources of infection on general graphs with $n$ vertices. External infections and recoveries occur at fixed rates $a,\kappa>0$. The per-edge infection is at rate $\lambda_n\ge0$. We give explicit upper and lower bounds for the mixing time starting from the worst possible initial condition. If the maximum total transmission infection rate stays below the sum of the external infection and recovery rates uniformly, we obtain $t_{\mathrm{mix}}^{(n)}(\epsilon)=\Theta(\log n)$. Under the stronger condition $\lambda_n^2\sum_x\deg(x)^2\to0$, we prove the sharp mixing time $t_{\mathrm{mix}}^{(n)}(\epsilon)=\frac{\log n}{2(a+\kappa)}+O_{a,\kappa,\epsilon}(1)$ and cutoff with a window of constant order. We apply these results to random graphs, including Erd\H{o}s--R\'enyi graphs, sparse random geometric graphs and Galton--Watson trees with binomial or Poisson offspring.
\end{abstract}

\maketitle

\section{Introduction}
\label{sec:intro}
In the standard susceptible--infected--susceptible (SIS) model, susceptible individuals become infected through interactions with infected individuals, and infected individuals become susceptible after recovery. In this paper, we study the SIS model on general finite graphs. Each configuration gives us the immunological statuses (susceptible or infected) of the vertices. We will consider a Continuous-Time Markov Chain (CTMC) to describe the dynamics.  Without external infection, the configuration in which every vertex is susceptible is an absorbing state. 
To make this state non-absorbing, we consider external sources of infection, so that a susceptible individual can become infected even when none of its neighbours is infected. The external source of infection could be environmental or a community source \cite{Sharker2023SIM}. We refer to this process as the noisy SIS model. The dynamics is described by three parts: external infections occur at rate $a>0$, infections transmitted along each edge occur at rate $\lambda_n\ge0$ ($n$ denoting the number of vertices in the graph), and infected individuals recover at rate $\kappa>0$. Under such dynamics, the chain is irreducible. We study quantitatively how quickly its configuration law approaches equilibrium.

The mixing time measures how long the Markov chain takes to approach its invariant measure in total variation distance, starting from the worst possible initial state. The cutoff phenomenon \cite{Aldous_Diaconis_cutoff,Diaconis1996cutoff} describes an abrupt transition when the distance from equilibrium drops from near one to near zero. The cutoff window is a time window, typically of much smaller order than the mixing time, when this abrupt drop occurs. For mixing time and the cutoff phenomenon, classical examples and general background material can be found in \cite{Aldous_Diaconis_cutoff,Diaconis1996cutoff,LP}. 
In a recent work \cite{HLR}, He, Luczak and Ross proved the cutoff phenomenon with a window of constant order for the infected-count process of the logistic SIS model with self-infection on the complete graph. 
Ramadas~\cite{Ramadas2014}
proved logarithmic mixing for the noisy voter model on arbitrary
finite graphs. Cox, Peres and Steif~\cite{CPS2016} established cutoff
with a window of constant order for noisy voter model with strong noise (constant order). Gon\c{c}alves, Jara, Marinho and Menezes~\cite{GJMM2026} obtained a 
cutoff profile for SSEP with reservoirs using
Yau's relative entropy method. Their approach is also relevant
to our cutoff proof.
Unlike the mean-field setting on complete graphs, the infection rates at different vertices depend on the locations of other infected vertices. We study the mixing time and cutoff phenomenon for the full configuration process in a more general framework. The graphs in this paper are allowed to be disconnected. The graphs can be dense, and their maximum degree may grow with $n$. 

We first establish upper and lower bounds for the mixing time on general graphs. Let $\Delta_n$ denote the maximum degree. In Theorem~\ref{thm:bounds}, we give explicit upper and lower bounds when $a+\kappa-\lambda_n\Delta_n>0$. In addition, if we require $\liminf_{n\to\infty}(a+\kappa-\lambda_n\Delta_n)>0$, then Corollary~\ref{cor:log} gives the mixing time has the exact order $\log n$.

We then prove a sharp mixing time and a cutoff phenomenon, under a stronger condition on the per-edge infection rates. If $\lambda_n^2\sum_{x\in V_n}\deg(x)^2\to0$, then Theorem~\ref{thm:cutoff-degree-squares} gives a sharp mixing time $\tmix(\epsilon)=\frac{\log n}{2(a+\kappa)}+O_{a,\kappa,\epsilon}(1)$
for every fixed $\epsilon\in(0,1)$.
In addition, the chain exhibits a cutoff with a window of constant order. The proof uses Yau's relative entropy method  \cite{JaraMenezes,XK2026}.

Finally, we apply these results to three types of random graphs, including Erd\H{o}s--R\'enyi graphs, random geometric graphs and Galton--Watson trees with binomial or Poisson offspring. Interestingly, bounding the maximum degree gives the mixing time of order $\log n$ with high probability. Under the additional infection rate conditions stated in Theorems~\ref{thm:mixing_erdos}, \ref{thm:mixing_geometric} and~\ref{thm:mixing_galton}, bounds on the sum of squared degrees give a cutoff with high probability. In each case, the leading cutoff time is $\log n/[2(a+\kappa)]$.

The paper is structured as follows. Section~\ref{sec:model} defines the model and states the mixing bounds. Sections~\ref{sec:upper} and~\ref{sec:lower} prove the upper and the lower bounds. Section~\ref{sec:cutoff} proves the sharp mixing time and cutoff phenomenon. Section~\ref{sec:random} is devoted to studying the mixing properties on random graph examples.

\section{The model and the mixing time}
\label{sec:model}

\subsection{Definition and dynamics of the model}
Let $n\in\mathbb N$ be the scaling parameter. We consider a finite undirected graph $G_n=(V_n,E_n)$ with $|V_n|=n$, where $V_n$ and $E_n$ denote the vertex set and the edge set, respectively. Write $x\sim y$ when $\{x,y\}\in E_n$. We denote by $\Delta_n$
the maximal degree of $G_n$. 
Each vertex $x\in V_n$ represents an individual at site $x$.
In a configuration $\sigma\in\Omega_n\defeq
\{0,1\}^{V_n}$, the individual at $x$ is \emph{infected} if $\sigma_x=1$ and
\emph{susceptible} if $\sigma_x=0$. The number of infected neighbours of $x$
and its infection rate are given by
\[
 n_{I,\sigma}(x)=\sum_{y\sim x}\sigma_y,\qquad
 p(\sigma,x)=a+\lambda_n n_{I,\sigma}(x).
\]
Here, the parameter $a>0$ is the external infection rate and $\lambda_n\ge0$ is the
infection rate per infected neighbour. $\lambda_n$ may depend on $n$. An infected individual recovers at rate $\kappa>0$. 
We define $\sigma^x$ as the configuration obtained by flipping the state at
$x$. We define the noisy SIS model on $G_n$ to be the continuous-time
Markov chain $(\sigma(t))_{t\ge0}$ with generator
\begin{equation}
 \label{eq:generator}
 (\cL_n f)(\sigma)=\sum_{x\in V_n}
   \bigl[(1-\sigma_x)p(\sigma,x)+\sigma_x\kappa\bigr]
   \bigl[f(\sigma^x)-f(\sigma)\bigr].
\end{equation}
For $f,g:\Omega_n\to\mathbb R$, we define the \emph{carr\'e du champ} operator associated with $\cL_n$ by $\Gamma_n(f,g)\defeq\cL_n(fg)-f\cL_n g-g\cL_n f$. Writing $\Gamma_n g\defeq\Gamma_n(g,g)$, we obtain
\begin{equation}
 \label{eq:carre-du-champ}
 (\Gamma_n g)(\sigma)=\sum_{x\in V_n}\bigl[(1-\sigma_x)p(\sigma,x)+\sigma_x\kappa\bigr]\bigl[g(\sigma^x)-g(\sigma)\bigr]^2.
\end{equation}
Because of the external infections, the chain is irreducible and has
a unique invariant measure, which we denote by $\mu_n$.

\subsection{Distance to stationarity}
For $t\ge0$ and $\eta\in\Omega_n$, denote the law of the Markov chain at
time $t$ starting from $\eta$ by $\mu_{\eta,t}^{(n)}$. The associated semigroup is denoted by 
$S_tf(\eta)=\E_\eta[f(\sigma(t))]$. We write $\nu(f)$ for the expectation
of $f$ under $\nu$ and $\Var_\nu(f)=\nu(f^2)-\nu(f)^2$ for its variance.
Given two probability measures $\mu$ and $\nu$ on $\Omega_n$, the
\emph{total variation distance} between them is
\[
 \TV{\nu}{\mu}=\frac12\sum_{\sigma\in\Omega_n}|\nu(\sigma)-\mu(\sigma)|
              =\max_{A\subseteq\Omega_n}|\nu(A)-\mu(A)|.
\]
For $\epsilon\in(0,1)$, the worst-case distance to stationarity and the
\emph{$\epsilon$-mixing time} are given by
\begin{align}
 \label{eq:mixing}
 \mathrm d^{(n)}(t)
   &\defeq \max_{\eta\in\Omega_n}\TV{\mu_{\eta,t}^{(n)}}{\mu_n},
 &
 \tmix(\epsilon)
   &\defeq\inf\{t\ge 0:\mathrm d^{(n)}(t)\le\epsilon\}.
\end{align}
We say that the family of Markov chains exhibits \emph{total variation cutoff} if 
\[
\tmix(\epsilon)/\tmix(1-\epsilon)\to1,
\]
as $n\to\infty$ for every fixed $\epsilon\in(0,1/2)$. The cutoff occurs at $t_n\to\infty$ if $\tmix(\epsilon)/t_n\to1$ for every fixed $\epsilon\in(0,1)$.
A sequence $w_n>0$ with $w_n=o(t_n)$ is a \emph{cutoff window} around $t_n$ if, for every fixed $\epsilon\in(0,1)$,
 $\tmix(\epsilon)=t_n+O_\epsilon(w_n)$.
A window of constant order corresponds to taking $w_n=1$.

\subsection{Mixing bounds}
Let $[u]_+=\max\{u,0\}$.
To state the bounds, we define
\begin{equation}
 \label{eq:rates}
 c_n=a+\kappa-\lambda_n\Delta_n,\qquad
 b_n=a+\kappa+\lambda_n\Delta_n,\qquad
 R_n=\max\{a+\lambda_n\Delta_n,\kappa\}.
\end{equation}

\begin{theorem}
 \label{thm:bounds}
Suppose $c_n>0$. Then, for every $\epsilon\in(0,1)$,
 \begin{equation}
  \label{eq:time-bounds}
  \frac{1}{2b_n}
     \left[\log\!\left(\frac{(1-\epsilon)c_n n}{16R_n}\right)\right]_+
  \ \le\ \tmix(\epsilon)
  \ \le\ \frac{\log n+\log(1/\epsilon)}{c_n}.
 \end{equation}
\end{theorem}
The proof of Theorem~\ref{thm:bounds} is split into two parts. The upper bound is proved  in Section~\ref{sec:upper}. The proof of the lower bound is given in Section~\ref{sec:lower}. 

\begin{corollary}
\label{cor:log}
Suppose that $\liminf_{n\to\infty}c_n>0$.
Then, for each fixed $\epsilon\in(0,1)$ and all sufficiently large $n$,
\begin{align}
\label{mixing_time_bound_corollary}
 \frac{\log n}{2b_n}-O(1)
 \le \tmix(\epsilon)
 \le \frac{\log n+\log(1/\epsilon)}{c_n},
\end{align}
where the constant in $O(1)$ depends only on $a,\kappa,\epsilon$
and $\liminf_{n\to\infty}c_n$.
Consequently, $\tmix(\epsilon)=\Theta(\log n)$.
\end{corollary}

\begin{proof}[Proof of Corollary~\ref{cor:log}]
Since $\liminf_{n\to\infty}c_n>0$, we have $1/c_n=O(1)$. For all
sufficiently large $n$, $0<c_n\le a+\kappa$,
$a+\kappa\le b_n\le 2(a+\kappa)$, and
$\max\{a,\kappa\}\le R_n\le 2a+\kappa$. Hence, for fixed $\epsilon$,
\[
 \log\!\left(\frac{(1-\epsilon)c_n n}{16R_n}\right)
 =\log n+O(1),
\]
where the constant depends only on $a,\kappa,\epsilon$ and
$\liminf_{n\to\infty}c_n$.
\end{proof}

\section{The upper bound}
\label{sec:upper}

Consider two copies of the chain started from configurations
$\eta,\xi\in\Omega_n$. Denote their configurations at time $t$ by
$X^\eta_t$ and $X^\xi_t$. We consider the partial order: if $\eta_x\le\xi_x$ for all $x\in V_n$, we
write $\eta\le\xi$. Let
\[
 F(\eta)=\sum_{x\in V_n}\eta_x,\qquad
 H(\eta,\xi)=\sum_{x\in V_n}|\xi_x-\eta_x|.
\]
The nonnegative integer $H(\eta,\xi)$ is the \emph{Hamming distance} between two configurations, which counts the number of vertices where
$\eta$ and $\xi$ disagree. We suppose that $\eta\le\xi$, then $H(\eta,\xi)=F(\xi)-F(\eta)$. 

The coupling is constructed as follows: All copies are constructed using the same independent Poisson clocks. At each vertex, a clock of rate $a$ sets its state to $1$, and a clock of rate $\kappa$ sets its state to $0$.
Along each edge $\{x,y\}\in E_n$, infection attempts occur at rate $\lambda_n$ in each direction. In each copy, an attempt from $y$ to $x$ infects $x$ if $y$ is infected and $x$ is susceptible. Otherwise, nothing changes.

We summarize the coupling in the following table.

\begin{center}
\renewcommand{\arraystretch}{1.2}
\begin{tabular}{ccl}
\hline
$(\eta_x,\xi_x)$ & Transition & Rate \\
\hline
$(0,0)$ & $(\eta,\xi)\to(\eta^x,\xi^x)$ & $p(\eta,x)$ \\
$(0,0)$ & $(\eta,\xi)\to(\eta,\xi^x)$   & $p(\xi,x)-p(\eta,x)$ \\
$(1,1)$ & $(\eta,\xi)\to(\eta^x,\xi^x)$ & $\kappa$ \\
$(0,1)$ & $(\eta,\xi)\to(\eta^x,\xi)$   & $p(\eta,x)$ \\
$(0,1)$ & $(\eta,\xi)\to(\eta,\xi^x)$   & $\kappa$ \\
\hline
\end{tabular}
\end{center}
Each transition preserves the coordinatewise order. Therefore, the coupling is monotone.

\begin{lemma}
 \label{lem:drift}
 For $\eta\le\xi$ and $t\ge 0$,
 \begin{equation}
  \label{eq:two-rates}
  e^{-b_n t}H(\eta,\xi)
  \ \le\ S_tF(\xi)-S_tF(\eta)
  =\E H(X^\eta_t,X^\xi_t)
  \ \le\ e^{-c_n t}H(\eta,\xi).
 \end{equation}
\end{lemma}

\begin{proof}[Proof of Lemma~\ref{lem:drift}]
 Let $e_{SI}(\sigma)$ be the number of edges with one susceptible
 and one infected endpoint. Then, the generator in \eqref{eq:generator}
 reduces to 
 \begin{equation}
  \label{eq:count-drift}
  (\cL_n F)(\sigma)=an-(a+\kappa)F(\sigma)+\lambda_n e_{SI}(\sigma).
 \end{equation}
 One flip at site $x$ changes the value of $e_{SI}(\sigma)$ by at most $\deg(x)$
 edges. Therefore, it follows that
 $|e_{SI}(\xi)-e_{SI}(\eta)|\le\Delta_n H(\eta,\xi)$.
 For $\eta\le\xi$, subtraction in \eqref{eq:count-drift} now gives
 \[
  -b_n H(\eta,\xi)\le (\cL_n F)(\xi)-(\cL_n F)(\eta)\le -c_n H(\eta,\xi).
 \]
 We define $h(t)=S_tF(\xi)-S_tF(\eta)$. It follows that
 $h(t)=\E H(X^\eta_t,X^\xi_t)$. Differentiating with respect to $t$,
 $h'(t)=\E[(\cL_n F)(X^\xi_t)-(\cL_n F)(X^\eta_t)]$.
 Therefore, we have that
 $-b_n h(t)\le h'(t)\le-c_n h(t)$, integrating which we get \eqref{eq:two-rates}.
\end{proof}

Let $\one$  and $\zero$  denote the all-infected and all-susceptible configurations, respectively. We are now ready to prove the upper bound in Theorem~\ref{thm:bounds}. 

\begin{proof}[Proof of the upper bound in Theorem~\ref{thm:bounds}]
Draw $\xi$ from $\mu_n$ independently of the
clocks. By stationarity, $X_t^\xi$ has law $\mu_n$ for every $t\ge0$. We run the four copies started from $\zero,\eta,\xi,\one$ using
the same clocks.  Monotonicity of this coupling keeps both $X_t^\eta$ and $X_t^\xi$
between $X_t^\zero$ and $X_t^\one$ coordinatewise. 
The coupling inequality \cite[Chapter~5]{LP} and
Lemma~\ref{lem:drift} yield the following estimate
\begin{align*}
 \TV{\mu_{\eta,t}^{(n)}}{\mu_n}
 &\le \PP(X_t^\eta\ne X_t^\xi)
 \le \E H(X_t^\eta,X_t^\xi)\\
 &\le \E H(X_t^\zero,X_t^\one)
 \le e^{-c_n t}H(\zero,\one)=ne^{-c_n t}.
\end{align*}
Maximising over $\eta$ then gives
$\mathrm d^{(n)}(t)\le\min\{1,ne^{-c_n t}\}$. Since $c_n>0$,
choosing $t=(\log n+\log(1/\epsilon))/c_n$ proves the upper bound in
\eqref{eq:time-bounds}.
\end{proof}

\section{The lower bound}
\label{sec:lower}

Take $m=\mu_n(F)$.
We suppose that the Markov chain starts from $\zero$. Couple it with a
copy $\xi$ with initial law $\mu_n$. Since
$\zero\le\xi$, the lower bound in Lemma~\ref{lem:drift} applies for
each $\xi$. Integrating over $\xi$ and using $\mu_n(S_tF)=m$ gives the
first inequality below: 
\[
  m-S_tF(\zero)\ge e^{-b_n t}m,\qquad
  S_tF(\one)-m\ge e^{-b_n t}(n-m).
\]
The second follows from
$\xi\le\one$.
Take $\eta_*=\zero$ if $m\ge n/2$, and $\eta_*=\one$ otherwise. We conclude that
\begin{equation}
 \label{eq:mean-gap}
 |S_tF(\eta_*)-m|\ge \tfrac{n}{2}e^{-b_n t}, 
 \qquad\text{for every }t\ge 0 .
\end{equation}

\begin{lemma}
 \label{lem:variance}
Suppose $c_n>0$. For any initial condition $\eta \in \Omega_n$ and any $t\ge 0$,
 \begin{equation}
  \label{eq:variances}
  \Var_{\mu_{\eta,t}^{(n)}}(F)\le \frac{nR_n}{2c_n}(1-e^{-2c_n t}),
  \qquad
  \Var_{\mu_n}(F)\le \frac{nR_n}{2c_n}.
 \end{equation}
\end{lemma}

\begin{proof}[Proof of Lemma~\ref{lem:variance}]
 For any $u \geq 0$, Lemma~\ref{lem:drift} gives that
 \begin{equation}
  \label{eq:gradient}
  |S_uF(\sigma^x)-S_uF(\sigma)|\le e^{-c_n u}.
 \end{equation}
 For $0\le s\le t$, we define $g_s=S_{t-s}F$. Since
 $\partial_s g_s=-\cL_n g_s$, differentiating $S_s(g_s^2)$ gives
 \[
  \frac{\mathrm d}{\mathrm ds} S_s(g_s^2)
    =S_s\!\left(\cL_n(g_s^2)-2g_s\cL_n g_s\right)=S_s(\Gamma_n g_s).
 \]
Therefore, by integrating from $0$ to $t$, we get the variance
 identity
 \[
  \Var_{\mu_{\eta,t}^{(n)}}(F)
    = S_t(F^2)(\eta)-(S_tF(\eta))^2
    = \int_0^t S_s(\Gamma_n g_s)(\eta)\,\mathrm ds .
 \]
By \eqref{eq:gradient} and taking expectation under $\mu_{\eta,s}^{(n)}$, we get
$S_s(\Gamma_n g_s)(\eta)\le nR_n e^{-2c_n(t-s)}$. Integrating over
$s\in[0,t]$ gives the first estimate in \eqref{eq:variances}.
For a fixed $n$, $\mu_{\eta,t}^{(n)}$ converges to $\mu_n$ as $t\to\infty$.
Since both $F$ and $F^2$ are bounded, passing to the limit gives the
second estimate, which concludes the proof. 
\end{proof}

We are now ready to present the proof of the lower bound in Theorem~\ref{thm:bounds}. 

\begin{proof}[Proof of the lower bound in Theorem~\ref{thm:bounds}]
 We define $m_t=\mu_{\eta_*,t}^{(n)}(F)$ and take
 $r=(m_t+m)/2$. If $m_t>m$, take $A_t=\{\sigma:F(\sigma)>r\}$.
 If $m_t<m$, take $A_t=\{\sigma:F(\sigma)<r\}$. By \eqref{eq:mean-gap}, if a configuration $\sigma$ sampled from
$\mu_{\eta_*,t}^{(n)}$ lies in $A_t^c$, then $F(\sigma)\le r<m_t$
when $m_t>m$, and $m_t<r\le F(\sigma)$ when $m_t<m$. In either case,
$|F(\sigma)-m_t|\ge |r-m_t|=|m_t-m|/2$.

Similarly, if a configuration $\sigma$ drawn from $\mu_n$ lies in $A_t$, then
$r$ lies strictly between $F(\sigma)$ and $m$, so
$|F(\sigma)-m|>|r-m|=|m_t-m|/2$.
 Applying Chebyshev's inequality to $F$,  
\begin{align*}
 \mu_{\eta_*,t}^{(n)}(A_t^c)
 &\le \mu_{\eta_*,t}^{(n)}\!\left(|F-m_t|\ge \frac{|m_t-m|}{2}\right)
  \le \frac{4\Var_{\mu_{\eta_*,t}^{(n)}}(F)}{|m_t-m|^2},\\
 \mu_n(A_t)
 &\le \mu_n\!\left(|F-m|\ge \frac{|m_t-m|}{2}\right)
  \le \frac{4\Var_{\mu_n}(F)}{|m_t-m|^2}.
\end{align*}
Therefore, it follows that
\begin{align*}
 \mathrm d^{(n)}(t)
 \geq 1-\mu_{\eta_*,t}^{(n)}(A_t^c)-\mu_n(A_t)
 \ge 1-\frac{4\bigl(\Var_{\mu_{\eta_*,t}^{(n)}}(F)
                    +\Var_{\mu_n}(F)\bigr)}{|m_t-m|^2}.
\end{align*}
 Inserting the estimates from \eqref{eq:mean-gap} and \eqref{eq:variances} yields
 \[
  \mathrm d^{(n)}(t)
   \ge 1-
      \frac{4[nR_n/(2c_n)](2-e^{-2c_n t})}{(n^2/4)e^{-2b_n t}}
   =1-\frac{8R_n}{c_n n}(2-e^{-2c_n t})e^{2b_n t}.
 \]
 Since
 $2-e^{-2c_n t}\le 2$, the distance is greater than $\epsilon$ at every
 $t\ge 0$ such that
 $
  t<\frac{1}{2b_n}\log\!\left(\frac{(1-\epsilon)c_n n}{16R_n}\right).
 $
 This proves the lower bound in \eqref{eq:time-bounds}.
\end{proof}

\section{Sharp mixing time and cutoff phenomenon}
\label{sec:cutoff}

We now prove the sharp mixing time and a cutoff phenomenon. The following theorem shows that the chain exhibits cutoff at $\log n/[2(a+\kappa)]$ with a window of constant order.

\begin{theorem}
\label{thm:cutoff-degree-squares}
Suppose that
\begin{equation}
 \label{eq:cutoff-degree-squares}
 \lambda_n^2\sum_{x\in V_n}\deg(x)^2\longrightarrow0.
\end{equation}
Then, for every fixed $\epsilon\in(0,1)$,
\begin{equation}
 \label{eq:cutoff-time}
 \tmix(\epsilon)=\frac{\log n}{2(a+\kappa)}+O_{a,\kappa,\epsilon}(1).
\end{equation}
The chain exhibits a cutoff at $\frac{\log n}{2(a+\kappa)}$ with a window of constant order.
\end{theorem}

\begin{proof}[Proof of Theorem~\ref{thm:cutoff-degree-squares}]
We let $\cL_{n,0}$ and $\mathrm d_0^{(n)}(t)$ denote the generator and worst-case distance for the process with the same $a,\kappa$, but $\lambda_n=0$. We first prove that
\begin{equation}
 \label{eq:cutoff-comparison}
 \sup_{t\ge0}\bigl|\mathrm d^{(n)}(t)-\mathrm d_0^{(n)}(t)\bigr|\le2\lambda_n\sqrt{\frac{\sum_{x\in V_n}\deg(x)^2}{a\min\{a,\kappa\}}}.
\end{equation}

We start both processes from the same product law $\bigotimes_{x\in V_n}\operatorname{Bern}(p_{x,0}^{(n)})$, with $p_{x,0}^{(n)}\in(0,1)$. Two processes have generators $\cL_n$ and $\cL_{n,0}$. We denote their time laws by $\mu_t^{(n)}$ and $\nu_t^{(n)}$, respectively. The process with generator $\cL_{n,0}$ has each site $x\in V_n$ evolving independently as $\lambda_n=0$.
This independent process remains a product law, with
\begin{align}
\label{reference_produc_laws}
 \nu_t^{(n)}=\bigotimes_{x\in V_n}\operatorname{Bern}(p_{x,t}^{(n)}),\qquad \partial_t p_{x,t}^{(n)}=a(1-p_{x,t}^{(n)})-\kappa p_{x,t}^{(n)}.
\end{align}
Define $f_t^{(n)}=\mu_t^{(n)}/\nu_t^{(n)}$ and $H_n(t)=\int_{\Omega_n}f_t^{(n)}\log f_t^{(n)}\,\mathrm d\nu_t^{(n)}$. Let $\cL_{n,t}^*$ be the adjoint of $\cL_n$ with respect to $\nu_t^{(n)}$. Applying Yau's relative entropy inequality \cite[Lemma~A.1]{JaraMenezes}, we get 
\begin{equation}
 \label{eq:Yau_entropy_upper_bound_F_n}
 H_n'(t)\le-\int_{\Omega_n}\Gamma_n\left(\sqrt{f_t^{(n)}}\right)\,\mathrm d\nu_t^{(n)}+\mathbb E_{\mu_t^{(n)}}\left[\cL_{n,t}^*\mathbf1-\partial_t\log\nu_t^{(n)}\right],
\end{equation}
where $\mathbf1$ denotes the constant function equal to one. For an application of Yau's relative entropy method, see \cite{XK2026}. Appendix~D of that paper contains computations of $\cL_{n,t}^*\mathbf1-\partial_t\log\nu_t^{(n)}$ as well. To compute $\cL_{n,t}^*\mathbf1-\partial_t\log\nu_t^{(n)}$, note that  $n_{I,\sigma^x}(x)=n_{I,\sigma}(x)$. Therefore, 
\[
 (\cL_{n,t}^*\mathbf1)(\sigma)=\sum_{x\in V_n}\frac{\sigma_x-p_{x,t}^{(n)}}{p_{x,t}^{(n)}(1-p_{x,t}^{(n)})}\left[(a+\lambda_n n_{I,\sigma}(x))(1-p_{x,t}^{(n)})-\kappa p_{x,t}^{(n)}\right].
\]
By \eqref{reference_produc_laws},  differentiating $\log\nu_t^{(n)}$ directly gives
\[
 \partial_t\log\nu_t^{(n)}(\sigma)=\sum_{x\in V_n}\frac{\sigma_x-p_{x,t}^{(n)}}{p_{x,t}^{(n)}(1-p_{x,t}^{(n)})}\left[a(1-p_{x,t}^{(n)})-\kappa p_{x,t}^{(n)}\right].
\]
Subtracting these identities, we obtain
\begin{equation}
 \label{eq:Yau_correction_explicit}
 \bigl(\cL_{n,t}^*\mathbf1-\partial_t\log\nu_t^{(n)}\bigr)(\sigma)=\lambda_n\sum_{x\in V_n}n_{I,\sigma}(x)\frac{\sigma_x-p_{x,t}^{(n)}}{p_{x,t}^{(n)}}.
\end{equation}
By \eqref{reference_produc_laws} and $f_t^{(n)}=\mu_t^{(n)}/\nu_t^{(n)}$, we get that
\[
 \mathbb E_{\mu_t^{(n)}}\left[\cL_{n,t}^*\mathbf1-\partial_t\log\nu_t^{(n)}\right]=\lambda_n\sum_{x\in V_n}\int_{\Omega_n}(1-\sigma_x)n_{I,\sigma}(x)\bigl[f_t^{(n)}(\sigma^x)-f_t^{(n)}(\sigma)\bigr]\,\mathrm d\nu_t^{(n)}.
\]
We denote by $\Gamma_{n,0}$  the carr\'e du champ operator associated with $\cL_{n,0}$, that is, 
\[
 \Gamma_{n,0}g(\sigma)=\sum_{x\in V_n}[a(1-\sigma_x)+\kappa\sigma_x]\bigl[g(\sigma^x)-g(\sigma)\bigr]^2.
\]
For any $t$ fixed, we write $g=\sqrt{f_t^{(n)}}$. Using $g(\sigma^x)^2-g(\sigma)^2=[g(\sigma^x)-g(\sigma)]^2+2g(\sigma)[g(\sigma^x)-g(\sigma)]$, we have
\[
 \begin{aligned}
 &\mathbb E_{\mu_t^{(n)}}\left[\cL_{n,t}^*\mathbf1-\partial_t\log\nu_t^{(n)}\right]\\
 &\quad=\int_{\Omega_n}(\Gamma_ng-\Gamma_{n,0}g)\,\mathrm d\nu_t^{(n)}\\
 &\qquad+2\lambda_n\sum_{x\in V_n}\int_{\Omega_n}(1-\sigma_x)n_{I,\sigma}(x)g(\sigma)\bigl[g(\sigma^x)-g(\sigma)\bigr]\,\mathrm d\nu_t^{(n)}.
 \end{aligned}
\]
Applying Young's inequality, we get 
\[
 2\lambda_n n_{I,\sigma}(x)g(\sigma)\bigl[g(\sigma^x)-g(\sigma)\bigr]\le\frac a2\bigl[g(\sigma^x)-g(\sigma)\bigr]^2+\frac{2\lambda_n^2}{a}n_{I,\sigma}(x)^2g(\sigma)^2.
\]
Since $n_{I,\sigma}(x)\le\deg(x)$ and $\int_{\Omega_n}g^2\,\mathrm d\nu_t^{(n)}=1$, it follows that
\begin{equation}
 \label{eq:Yau_correction_bound}
 \mathbb E_{\mu_t^{(n)}}\left[\cL_{n,t}^*\mathbf1-\partial_t\log\nu_t^{(n)}\right]\le\int_{\Omega_n}\Gamma_ng\,\mathrm d\nu_t^{(n)}-\frac12\int_{\Omega_n}\Gamma_{n,0}g\,\mathrm d\nu_t^{(n)}+\frac{2\lambda_n^2}{a}\sum_{x\in V_n}\deg(x)^2.
\end{equation}
From Yau's relative entropy estimate \eqref{eq:Yau_entropy_upper_bound_F_n}, we get 
\begin{equation}
 \label{eq:Yau_entropy_dissipation}
 H_n'(t)\le-\frac12\int_{\Omega_n}\Gamma_{n,0}g\,\mathrm d\nu_t^{(n)}+\frac{2\lambda_n^2}{a}\sum_{x\in V_n}\deg(x)^2.
\end{equation}

It remains to bound the term $\int_{\Omega_n}\Gamma_{n,0}g\,\mathrm d\nu_t^{(n)}$. Note that $\operatorname{Ent}_\pi(h^2)=\pi(h^2\log h^2)-\pi(h^2)\log\pi(h^2)$. In particular, $H_n(t)=\operatorname{Ent}_{\nu_t^{(n)}}(g^2)$. For $q\in(0,1)$ and $v,r>0$, we define
\[
 \Phi(r)=(1-q)v^2\log v^2+qr^2\log r^2-[(1-q)v^2+qr^2]\log[(1-q)v^2+qr^2].
\]
Note that $\Phi(v)=\Phi'(v)=0$.  Differentiating directly gives $ \Phi''(r) \leq 1$. Moreover, by 
Taylor's formula,  $\Phi(r)\le(r-v)^2/2$. Taking $v=h(0)$ and $r=h(1)$, and it follows that
$
 \operatorname{Ent}_{\operatorname{Bern}(q)}(h^2)\le\frac12[h(1)-h(0)]^2,\ h\ge0.
$
Since $\nu_t^{(n)}$ is a product law, by \cite[Proposition~5.6]{Ledoux2001} we conclude 
\[
 \operatorname{Ent}_{\nu_t^{(n)}}(h^2)\le\frac12\sum_{x\in V_n}\int_{\Omega_n}[h(\sigma^x)-h(\sigma)]^2\,\mathrm d\nu_t^{(n)}.
\]
At each vertex $x$, the independent process has transition rate satisfying $a(1-\sigma_x)+\kappa\sigma_x\ge\min\{a,\kappa\}$.
Therefore,
\begin{equation*}
 \operatorname{Ent}_{\nu_t^{(n)}}(h^2)\le\frac{1}{2\min\{a,\kappa\}}\int_{\Omega_n}\Gamma_{n,0}h\,\mathrm d\nu_t^{(n)}.
\end{equation*}
Recall the fact that $H_n(t)=\operatorname{Ent}_{\nu_t^{(n)}}(g^2)$. Applying the above inequality to $g$, gives
\begin{align*}
 H_n'(t)\le-\min\{a,\kappa\}H_n(t)+\frac{2\lambda_n^2}{a}\sum_{x\in V_n}\deg(x)^2.
\end{align*}
The two processes have the same initial law, so $H_n(0)=0$. Integration gives
\[
 H_n(t)\le\frac{2\lambda_n^2\sum_{x\in V_n}\deg(x)^2}{a\min\{a,\kappa\}}\left(1-e^{-\min\{a,\kappa\}t}\right).
\]
Applying Pinsker's inequality, we get 
\[
 \TV{\mu_t^{(n)}}{\nu_t^{(n)}}\le\sqrt{\frac12H_n(t)}\le\lambda_n\sqrt{\frac{\sum_{x\in V_n}\deg(x)^2}{a\min\{a,\kappa\}}}.
\]

Let $\nu_{\eta,t}^{(n)}$ and $\nu_n$ denote the law from $\eta$ at time $t$ and the invariant law of the independent process. Set $p_{x,0}^{(n)}=r+(1-2r)\eta_x$. By letting $r\downarrow0$, we can extend the estimate to every deterministic initial configuration in total variation distance. For each fixed $n$, sending $t\to\infty$ gives the same bound between two invariant measures $\mu_n$ and $\nu_n$. Therefore, to prove \eqref{eq:cutoff-comparison}, the following  is enough, which holds by the triangle inequality, 
\[
 \bigl|\mathrm d^{(n)}(t)-\mathrm d_0^{(n)}(t)\bigr|\le\max_{\eta\in\Omega_n}\TV{\mu_{\eta,t}^{(n)}}{\nu_{\eta,t}^{(n)}}+\TV{\mu_n}{\nu_n}\le2\lambda_n\sqrt{\frac{\sum_{x\in V_n}\deg(x)^2}{a\min\{a,\kappa\}}}.  
\]

When $\lambda_n=0$, each vertex is updated independently with rate $a+\kappa$ from $\operatorname{Bern}(a/(a+\kappa))$. The single-site update kernel is reversible and has eigenvalues $1$ and $0$, while its continuous-time generator has eigenvalues $0$ and $-(a+\kappa)$. This allows us to take $\gamma=1$ in \cite[Theorem~20.7 and Corollary~20.8]{LP}. The independent process has total updating rate $n(a+\kappa)$. Therefore, by rescaling time, we get 
\[
 \inf\{t\ge0:\mathrm d_0^{(n)}(t)\le\epsilon\}=\frac{\log n}{2(a+\kappa)}+O_{a,\kappa,\epsilon}(1),\qquad \epsilon\in(0,1).
\]
Under the condition in \eqref{eq:cutoff-degree-squares}, the error in \eqref{eq:cutoff-comparison} tends to zero, which proves
\eqref{eq:cutoff-time}. 
Finally, cutoff follows from \cite[Corollary~20.8]{LP}, \eqref{eq:cutoff-time}, and \eqref{eq:cutoff-comparison}. This concludes the proof. 
\end{proof}

\section{Mixing properties and cutoff on random graphs}
\label{sec:random}
In this section, we study the mixing time and cutoff of the noisy SIS model on random graphs. Write $\PP_G$ for the law of the sampled graphs. By cutoff at a location $t_0(n)\to \infty$ with a window of constant order with high probability, we mean that, for every fixed $\epsilon\in(0,1)$, there exists some deterministic $C_\epsilon<\infty$ satisfying
\[
 \PP_G\!\left(\left|\tmix(\epsilon)-t_0(n)\right|\le C_\epsilon\right)\longrightarrow1.
\]
In particular, $\tmix(\epsilon)/\tmix(1-\epsilon)\xrightarrow{\PP_G}1$ for every fixed $\epsilon\in(0,1/2)$.
\subsection{Mixing time on Erd\H{o}s--R\'enyi graphs}
\label{subsec:erdos-renyi}

Consider the Erd\H{o}s--R\'enyi graph $G_n\sim G(n,p_n)$. Each edge is included in the graph with probability $p_n$ independently. We write $d_n$ as the expected degree of a vertex, and we have $d_n=(n-1)p_n$ . We first sample
$G_n$ and then run the continuous-time noisy SIS chain on this sampled graph.
Let $\mu_n$ and $\tmix(\epsilon)$ refer to the invariant
measure and mixing time conditional on the realised graph.

\begin{theorem}
\label{thm:mixing_erdos}
Suppose that $d_n/\log n\to\infty$ and that the per-edge infection rates $\lambda_n\ge0$ satisfy $L\defeq\limsup_{n\to\infty}\lambda_n d_n<a+\kappa$. Then, for every fixed $\epsilon\in(0,1)$, there exists a constant $C=C(a,\kappa,L,\epsilon)$ such that, with $\PP_G$-probability $1-o(1)$, $c_n\ge(a+\kappa-L)/2$ and
\[
 \frac{\log n}{2b_n}-C\le\tmix(\epsilon)\le\frac{\log n+\log(1/\epsilon)}{c_n}.
\]
In particular, $\tmix(\epsilon)=\Theta(\log n)$ with high probability.
Under the additional condition $\sqrt n\lambda_n d_n\to0$, with high probability, the chain exhibits cutoff at $\log n/[2(a+\kappa)]$ with a window of constant order.
\end{theorem}

\begin{proof}[Proof of Theorem~\ref{thm:mixing_erdos}]
For each vertex $x$, $\deg(x)$ has a $\mathrm{Bin}(n-1,p_n)$
distribution with mean $d_n=(n-1)p_n$. We take
$\rho_n=\sqrt{6\log n/d_n}$. Since $d_n/\log n\to\infty$, we have
$0<\rho_n\le1$, for sufficiently large $n$. By a Chernoff bound and a union bound over the $n$ vertices, we get 
\[
 \PP_G\bigl(\Delta_n>(1+\rho_n)d_n\bigr)
 \le n\exp\!\left(-\frac{\rho_n^2 d_n}{3}\right)
 =\frac{1}{n}.
\]

Since $\rho_n\to0$ and $\limsup_{n\to\infty}\lambda_n d_n=L<a+\kappa$,
we have that $(1+\rho_n)\lambda_n d_n\le(a+\kappa+L)/2$
for all $n$ large enough.
We fix a realisation with $\Delta_n\le(1+\rho_n)d_n$. On this event, we obtain
\[
 c_n\ge a+\kappa-(1+\rho_n)\lambda_n d_n
 \ge \frac{a+\kappa-L}{2}>0.
\]
The upper bound follows from Theorem~\ref{thm:bounds}. For the
lower bound, \eqref{eq:time-bounds} gives
\[
 \tmix(\epsilon)
 \ge \frac{\log n}{2b_n}
     -\frac{1}{2b_n}\log\!\left(\frac{16R_n}{(1-\epsilon)c_n}\right).
\]
On the same event $\Delta_n\le(1+\rho_n)d_n$, we have
$a+\kappa\le b_n\le2(a+\kappa)$ and $R_n\le2a+\kappa$.
Together with the lower bound on $c_n$, these inequalities allow us to take
 $C=\frac{1}{2(a+\kappa)}
   \log\!\left(\frac{32(2a+\kappa)}
                     {(1-\epsilon)(a+\kappa-L)}\right)$ in \eqref{mixing_time_bound_corollary}.
The two bounds give $\tmix(\epsilon)=\Theta(\log n)$.

Suppose now that $\sqrt n\lambda_n d_n\to0$. On the same event $\Delta_n\le(1+\rho_n)d_n$, we have
\[
 \lambda_n^2\sum_{x\in V_n}\deg(x)^2\le n\lambda_n^2\Delta_n^2\le(1+\rho_n)^2n\lambda_n^2d_n^2\longrightarrow0.
\]
By Theorem~\ref{thm:cutoff-degree-squares}, \eqref{eq:cutoff-time} holds, and with $\PP_G$-probability $1-o(1)$, the chain exhibits cutoff.
\end{proof}

\subsection{Mixing time on random geometric graphs}
\label{subsec:geometric}
Let $G_n$ be the graph obtained by sampling $n$ points independently and uniformly on $[0,1]^2$ and connecting two points if their Euclidean distance is at most $r_n$. Suppose that $n\pi r_n^2\to\tau\in(0,\infty)$.
The maximum degree satisfies
 $\frac{\Delta_n\log\log n}{\log n}
 \xrightarrow{\PP_G}1$
by \cite[Theorem~6.10]{PenroseRGG}. 
Recent work on random geometric graphs includes estimates for the
mixing time of random walks on the giant component
\cite{HSRGG,KMMRGG}. In \cite{HMSRGG}, the
distinguishability of spherical random geometric graphs from Erd\H{o}s--R\'enyi graphs is studied.
Here, we consider the mixing time of the noisy SIS chain on the full
graph, conditioned on the sampled points.

\begin{theorem}
\label{thm:mixing_geometric}
Suppose that $n\pi r_n^2\to\tau\in(0,\infty)$ and that the per-edge infection rates $\lambda_n\ge0$ satisfy
 $L\defeq\limsup_{n\to\infty}\lambda_n\frac{\log n}{\log\log n}<a+\kappa$.
For every fixed $\epsilon\in(0,1)$, there exists a constant $C=C(a,\kappa,L,\epsilon)$ such that, with $\PP_G$-probability $1-o(1)$, $c_n\ge(a+\kappa-L)/2$ and
\[
 \frac{\log n}{2b_n}-C\le\tmix(\epsilon)\le\frac{\log n+\log(1/\epsilon)}{c_n}.
\]
In particular, $\tmix(\epsilon)=\Theta(\log n)$ with high probability. Under the additional condition $\sqrt n\lambda_n \to0$, with high probability, the chain exhibits cutoff at $\log n/[2(a+\kappa)]$ with a window of constant order.
\end{theorem}

\begin{proof}[Proof of Theorem~\ref{thm:mixing_geometric}]
Since $nr_n^2\to\tau/\pi\in(0,\infty)$, by Theorem~6.10
of \cite{PenroseRGG}, we have
 $\frac{\Delta_n\log\log n}{\log n}
 \xrightarrow{\PP_G}1$. Now we choose $\rho>0$ such that \[
 (1+\rho)L<(a+\kappa+L)/2.
\] With probability $1-o(1)$, $\Delta_n\le(1+\rho)\log n/\log\log n$. It follows that for all sufficiently large $n$,
 $\lambda_n\Delta_n\le(1+\rho)\lambda_n\frac{\log n}{\log\log n}\le\frac{a+\kappa+L}{2}$.
Consequently, $c_n\ge(a+\kappa-L)/2>0$, $a+\kappa\le b_n\le2(a+\kappa)$ and $R_n\le2a+\kappa$.

The upper bound follows directly from Theorem~\ref{thm:bounds}.
Substituting these bounds into \eqref{eq:time-bounds} gives the
stated inequalities with
$C=\frac{1}{2(a+\kappa)}
  \log\!\left(\frac{32(2a+\kappa)}{(1-\epsilon)(a+\kappa-L)}\right)$ in \eqref{mixing_time_bound_corollary}.
Therefore, it follows that $\tmix(\epsilon)=\Theta(\log n)$.

Suppose now that $\sqrt n\lambda_n\to0$. We condition on the location $u$ of a vertex, and the remaining $n-1$ points are independent and uniform in $[0,1]^2$. Each point is adjacent to this vertex with probability $\operatorname{Area}(B(u,r_n)\cap[0,1]^2)\le\pi r_n^2$, where $B(u,r_n)$ is the Euclidean disc of radius $r_n$ centred at $u$. It follows that the degree is binomial with conditional mean at most $(n-1)\pi r_n^2$. Since $n\pi r_n^2\to\tau$, we obtain
 $\E_G\sum_{x\in V_n}\deg(x)^2\le n\left[(n-1)\pi r_n^2+\bigl((n-1)\pi r_n^2\bigr)^2\right]=O(n)$.
Thus,  $\E_G[\lambda_n^2\sum_{x\in V_n}\deg(x)^2]=O(n\lambda_n^2)\to0$, and by Markov's inequality, we conclude
 \[
 \lambda_n^2\sum_{x\in V_n}\deg(x)^2\xrightarrow{\PP_G}0.
\]
By Theorem~\ref{thm:cutoff-degree-squares}, \eqref{eq:cutoff-time} holds, and with $\PP_G$-probability $1-o(1)$, the chain exhibits cutoff.
\end{proof}

\subsection{Mixing property on Galton--Watson trees}
\label{subsec: galton watson}
We start a Galton--Watson tree $T$ from a single individual. Each individual independently produces offspring with the same law. We
consider two different laws: the binomial and Poisson cases. In the binomial case, offspring are with law
$\mathrm{Bin}(n,p_n)$. In the Poisson case, the offspring distribution is
$\mathrm{Poisson}(\theta)$ with fixed mean $\theta>1$. 
For sufficiently large values of $n$, we condition on the event $\{|T|=\infty\}$ that the tree survives. Let $T_n$ be the subtree induced by the first $n$ vertices of $T$ in breadth-first order. We denote its conditional probability law by $\PP_G$. Given a realisation of $T_n$, we run the
continuous-time noisy SIS chain on this tree and write $\tmix(\epsilon)$
for its mixing time.

\begin{theorem}
\label{thm:mixing_galton}
In the binomial case, take $d_n=np_n$ satisfying $d_n/\log n\to\infty$. In the Poisson case, take $d_n=\log n/\log\log n$. We suppose that the per-edge infection rates $\lambda_n\ge0$ satisfy
 $L\defeq\limsup_{n\to\infty}\lambda_n d_n<a+\kappa$.
Then, for every fixed $\epsilon\in(0,1)$, there exists a constant $C=C(a,\kappa,L,\epsilon)$ such that, with $\PP_G$-probability $1-o(1)$, $c_n\ge(a+\kappa-L)/2$ and
\[
 \frac{\log n}{2b_n}-C\le\tmix(\epsilon)\le\frac{\log n+\log(1/\epsilon)}{c_n}.
\]
In particular, $\tmix(\epsilon)=\Theta(\log n)$ with high probability. Under the additional condition $\lambda_n\sqrt{nd_n}\to0$ in the binomial case, or $\sqrt n\lambda_n\to0$ in the Poisson case, with high probability, the chain exhibits cutoff at $\log n/[2(a+\kappa)]$ with a window of constant order.
\end{theorem}

\begin{proof}[Proof of Theorem~\ref{thm:mixing_galton}]
Before conditioning, we take an independent sequence $X_1,\allowbreak X_2,\allowbreak\ldots$
with the given offspring law. Under breadth-first order,
we use $X_i$ as the number of children of the $i$-th vertex.
In the tree $T_n$, this vertex has at most $X_i$ children and at most one
parent. Therefore, whenever the tree $T$ has at least $n$ vertices,
 $\Delta_n\le 1+\max_{1\le i\le n}X_i$.
Let us now choose a  small enough $\rho\in(0,1)$ such that
\[
 (1+\rho)L<(a+\kappa+L)/2.
\]

\emph{Binomial offspring.}
We define $q_n=\PP(|T|<\infty)$ as the probability of extinction. The number $q_n$ is the smallest fixed point of $f_n(s)=(1-p_n+p_ns)^n$ in $[0,1]$. Since $d_n\to\infty$, we have that $f_n(1/2)\le e^{-d_n/2}<1/2$, for $n$ large enough. Iteration from zero therefore gives $q_n=f_n(q_n)\le f_n(1/2)\le e^{-d_n/2}$. Hence, $\PP(|T|=\infty)=1-o(1)$. Using a Chernoff bound and a union bound before conditioning, we get
\[
 \PP_G\bigl(\Delta_n>1+(1+\rho)d_n\bigr)
 \le\frac{n\exp(-\rho^2d_n/3)}{\PP(|T|=\infty)}=o(1).
\]

\emph{Poisson offspring.}
We take $k_n=(1+\rho)d_n$. For all sufficiently large $n$, $k_n>\theta$,
and the Poisson Chernoff bound gives
$\PP(X_1>k_n)\le(e\theta/k_n)^{k_n}$. Since
$k_n\log(k_n/(e\theta))=(1+\rho+o(1))\log n$, we obtain
\[
 n\PP(X_1>k_n)\le n^{-\rho+o(1)}=o(1).
\]
The survival probability $\PP(|T|=\infty)$ is positive and independent
of $n$, since we have $\theta>1$. See \cite[Chapter~3]{Hofstad2017RGCNv1}.
Therefore, by the union bound, we conclude 
\[
 \PP_G\bigl(\Delta_n>1+(1+\rho)d_n\bigr)
 \le \frac{n\PP(X_1>k_n)}{\PP(|T|=\infty)}=o(1).
\]

In either case, fix a realisation with $\Delta_n\le 1+(1+\rho)d_n$.
Since $d_n\to\infty$ and $\lambda_n d_n$ is bounded, we have
$\lambda_n\to0$. By the choice of $\rho$, for all sufficiently large
$n$,
\[
 \lambda_n\Delta_n
 \le \lambda_n+(1+\rho)\lambda_n d_n
 \le \frac{a+\kappa+L}{2}.
\]
Hence $c_n\ge(a+\kappa-L)/2$, $a+\kappa\le b_n\le2(a+\kappa)$,
and $R_n\le2a+\kappa$. Substituting these bounds into
\eqref{eq:time-bounds} gives the stated inequalities with
 $C=\frac{1}{2(a+\kappa)}
   \log\!\left(\frac{32(2a+\kappa)}
                     {(1-\epsilon)(a+\kappa-L)}\right)$ in 
                     \eqref{mixing_time_bound_corollary}.
It follows that
$\tmix(\epsilon)=\Theta(\log n)$.

We now prove cutoff under the additional conditions. In the binomial case, suppose that $\lambda_n\sqrt{nd_n}\to0$. Since $T_n$ is a tree with $n$ vertices, $\sum_{x\in V_n}\deg(x)=2(n-1)$. On the same event $\Delta_n\le1+(1+\rho)d_n$, we have $\Delta_n\le3d_n$ for all sufficiently large $n$. Therefore,
\[
 \lambda_n^2\sum_{x\in V_n}\deg(x)^2\le2(n-1)\lambda_n^2\Delta_n\le6n\lambda_n^2d_n\longrightarrow0.
\]

In the Poisson case, suppose that $\sqrt n\lambda_n\to0$. On the survival event, the degree of the $i$th vertex in $T_n$ is at most $1+X_i$, so $\sum_{x\in V_n}\deg(x)^2\le\sum_{i=1}^n(1+X_i)^2$. Since $\E[(1+X_1)^2]=1+3\theta+\theta^2$, we obtain
\[
 \E_G\sum_{x\in V_n}\deg(x)^2
 \le\frac{\E\left[\sum_{i=1}^n(1+X_i)^2\right]}{\PP(|T|=\infty)}
 =\frac{n(1+3\theta+\theta^2)}{\PP(|T|=\infty)}
 =O(n).
\]
Thus, $\E_G[\lambda_n^2\sum_{x\in V_n}\deg(x)^2]=O(n\lambda_n^2)\to0$. Finally, by Markov's inequality, 
\[
 \lambda_n^2\sum_{x\in V_n}\deg(x)^2\xrightarrow{\PP_G}0.
\]
 In either case, by Theorem~\ref{thm:cutoff-degree-squares}, \eqref{eq:cutoff-time} holds, and with $\PP_G$-probability $1-o(1)$, the chain exhibits cutoff.
\end{proof}

\section*{Acknowledgements}
The authors used AI tools (GPT-6 Astra) to check grammar, and improve the writing. The mathematical results and proofs were developed independently by the authors.
This work was supported in part by the Engineering and Physical Sciences Research Council (EPSRC) of UK Research and Innovation (UKRI) [EP/Y027795/1] and the National Science Foundation (NSF) [DMS-2246815].

\nocite{KX}

\bibliographystyle{amsplain}
\bibliography{references}

@article{HLR,
  author  = {He, R. and Luczak, M. and Ross, N.},
  title   = {Cutoff for the logistic {SIS} epidemic model with self-infection},
  journal = {Adv. in Appl. Probab.},
  volume  = {57},
  year    = {2025},
  pages   = {1456--1483},
  doi     = {10.1017/apr.2025.16},
  note    = {\href{https://doi.org/10.1017/apr.2025.16}{doi:10.1017/apr.2025.16}}
}

@misc{KX,
  author        = {KhudaBukhsh, W. R. and Xiang, Y.},
  title         = {Mixing time for a noisy {SIS} model on graphs},
  year          = {2025},
  eprint        = {2501.07738v1},
  archivePrefix = {arXiv},
  url           = {https://arxiv.org/abs/2501.07738v1},
  howpublished  = {\ARXIV{2501.07738v1}}
}

@book{LP,
  author    = {Levin, D. A. and Peres, Y.},
  title     = {Markov Chains and Mixing Times},
  edition   = {Second},
  publisher = {American Mathematical Society},
  address   = {Providence, RI},
  year      = {2017},
  mrnumber  = {3726904},
  doi       = {10.1090/mbk/107},
  note      = {With contributions by E. L. Wilmer. \href{https://doi.org/10.1090/mbk/107}{doi:10.1090/mbk/107}}
}

@book{Hofstad2017RGCNv1,
  author    = {van der Hofstad, R.},
  title     = {Random Graphs and Complex Networks},
  volume    = {1},
  publisher = {Cambridge University Press},
  address   = {Cambridge},
  year      = {2017},
  mrnumber  = {3617364},
  doi       = {10.1017/9781316779422},
  note      = {\href{https://doi.org/10.1017/9781316779422}{doi:10.1017/9781316779422}}
}

@book{PenroseRGG,
  author    = {Penrose, M.},
  title     = {Random Geometric Graphs},
  publisher = {Oxford University Press},
  address   = {Oxford},
  year      = {2003},
  mrnumber  = {1986198},
  doi       = {10.1093/acprof:oso/9780198506263.001.0001},
  note      = {\href{https://doi.org/10.1093/acprof:oso/9780198506263.001.0001}{doi:10.1093/acprof:oso/9780198506263.001.0001}}
}

@misc{HSRGG,
  author        = {Haaland, M. H. and {\v{S}}arkovi{\'c}, A.},
  title         = {Mixing time of the random walk on the giant component of the random geometric graph},
  year          = {2025},
  eprint        = {2512.21322},
  archivePrefix = {arXiv},
  url           = {https://arxiv.org/abs/2512.21322},
  howpublished  = {\ARXIV{2512.21322}}
}

@misc{KMMRGG,
  author        = {Kiwi, M. and Martinez, C. and Mitsche, D.},
  title         = {Mixing time and isoperimetry in random geometric graphs},
  year          = {2025},
  eprint        = {2510.19951},
  archivePrefix = {arXiv},
  url           = {https://arxiv.org/abs/2510.19951},
  howpublished  = {\ARXIV{2510.19951}}
}

@misc{HMSRGG,
  author        = {Hunter, Z. and Milojevi{\'c}, A. and Sudakov, B.},
  title         = {Distinguishability threshold for random geometric graphs},
  year          = {2026},
  eprint        = {2607.22480},
  archivePrefix = {arXiv},
  url           = {https://arxiv.org/abs/2607.22480},
  howpublished  = {\ARXIV{2607.22480}}
}

@misc{JaraMenezes,
  author        = {Jara, M. and Menezes, O.},
  title         = {Non-equilibrium fluctuations of interacting particle systems},
  year          = {2018},
  eprint        = {1810.09526},
  archivePrefix = {arXiv},
  url           = {https://arxiv.org/abs/1810.09526},
  howpublished  = {\ARXIV{1810.09526}}
}

@book{Ledoux2001,
  author    = {Ledoux, M.},
  title     = {The Concentration of Measure Phenomenon},
  series    = {Mathematical Surveys and Monographs},
  volume    = {89},
  publisher = {American Mathematical Society},
  address   = {Providence, RI},
  year      = {2001},
  mrnumber  = {1849347},
  doi       = {10.1090/surv/089},
  note      = {\href{https://doi.org/10.1090/surv/089}{doi:10.1090/surv/089}}
}

@misc{XK2026,
  author        = {Xiang, Y. and KhudaBukhsh, W. R.},
  title         = {Cutoff profiles and boundary effects in the {Togashi--Kaneko} model of autocatalytic reactions with weak mutations},
  year          = {2026},
  eprint        = {2610.17775},
  archivePrefix = {arXiv},
  url           = {https://arxiv.org/abs/2610.17775},
  howpublished  = {\ARXIV{2610.17775}}
}

@article{Aldous_Diaconis_cutoff,
  author  = {Aldous, D. and Diaconis, P.},
  title   = {Shuffling cards and stopping times},
  journal = {Amer. Math. Monthly},
  volume  = {93},
  year    = {1986},
  pages   = {333--348},
  mrnumber = {0841111},
  doi     = {10.1080/00029890.1986.11971821},
  note    = {\href{https://doi.org/10.1080/00029890.1986.11971821}{doi:10.1080/00029890.1986.11971821}}
}

@article{Diaconis1996cutoff,
  author  = {Diaconis, P.},
  title   = {The cutoff phenomenon in finite {Markov} chains},
  journal = {Proc. Natl. Acad. Sci. USA},
  volume  = {93},
  year    = {1996},
  pages   = {1659--1664},
  mrnumber = {1374011},
  doi     = {10.1073/pnas.93.4.1659},
  note    = {\href{https://doi.org/10.1073/pnas.93.4.1659}{doi:10.1073/pnas.93.4.1659}}
}

@article{Sharker2023SIM,
  author  = {Sharker, Y. and Diallo, Z. and KhudaBukhsh, W. R. and Kenah, E.},
  title   = {{Pairwise Accelerated Failure Time Regression Models for Infectious Disease Transmission in Close-Contact Groups With External Sources of Infection}},
  journal = {Stat. Med.},
  volume  = {43},
  year    = {2024},
  pages   = {5138--5154},
  doi     = {10.1002/sim.10226},
  note    = {\href{https://doi.org/10.1002/sim.10226}{doi:10.1002/sim.10226}}
}

@article{Ramadas2014,
  author  = {Ramadas, Harishchandra},
  title   = {Mixing of the noisy voter model},
  journal = {Electronic Communications in Probability},
  year    = {2014},
  volume  = {19},
  number  = {17},
  pages   = {1--9},
  doi     = {10.1214/ECP.v19-2968}
}

@article{CPS2016,
  author  = {Cox, J. Theodore and Peres, Yuval and Steif, Jeffrey E.},
  title   = {Cutoff for the noisy voter model},
  journal = {The Annals of Applied Probability},
  year    = {2016},
  volume  = {26},
  number  = {2},
  pages   = {917--932},
  doi     = {10.1214/15-AAP1108}
}

@article{GJMM2026,
  author  = {Gon{\c{c}}alves, Patr{\'i}cia and Jara, Milton
             and Marinho, Rodrigo and Menezes, Ot{\'a}vio},
  title   = {Sharp convergence to equilibrium for the {SSEP}
             with reservoirs},
  journal = {Annales de l'Institut Henri Poincar{\'e},
             Probabilit{\'e}s et Statistiques},
  year    = {2026},
  volume  = {62},
  number  = {1},
  pages   = {146--165},
  doi     = {10.1214/23-AIHP1419}
}

\end{document}